\documentclass[reqno,12pt,centertags]{amsart}
\usepackage{amsmath,amsthm,amscd,amssymb,latexsym,upref,stmaryrd,geometry,setspace}

\usepackage[CJKbookmarks]{hyperref}
\usepackage[perpage,symbol*]{footmisc}
\newcommand*{\mailto}[1]{\href{mailto:#1}{\nolinkurl{#1}}}

\newcommand{\beq}{\begin{equation}}
\newcommand{\eeq}{\end{equation}}
\newcommand{\ba}{\begin{align}}
\newcommand{\ea}{\end{align}}

\renewcommand{\Im}{\text{\rm Im}}

 \allowdisplaybreaks
\numberwithin{equation}{section}

\newtheorem{theorem}{Theorem}[section]
\newtheorem{lemma}[theorem]{Lemma}

\theoremstyle{definition}

\newtheorem{remark}{Remark}[section]
\makeatletter
\renewcommand{\section}{\@startsection{section}{1}{0mm}
  {-\baselineskip}{0.5\baselineskip}{\bf\leftline}}
\makeatother
\begin{document}


\noindent{\title[]
{\Large{Half Inverse Sturm-Liouville Problem on Three-Star Graph with a Discontinuity}}}

\author[]
{{{Kang Lv \ }} \\ \small{(Nanjing University of Science \& Technology ,  Nanjing, China )}}

\maketitle
\noindent {\textbf{\small{Abstract}}: \small{We consider the inverse Sturm-Liouville problem with one discontinuous point on three-star graph, deduced the distribution of the eigenvalues, and proved that one spectrum could uniquely determine the unknown potential and jump information when we know the whole potentials on two edges and half potential on another edge. }}


\footnotetext{lvkang201905@outlook.com}
\section{\large{Introduction}}
 O. Hald[1] considered the half inverse Sturm-Liouville problem with a discontinuous point, proved that one spectrum with half potential could uniquely decide the unknown potential and all jump information. C. Willis [8] studied the half inverse Sturm-Liouville problem with two discontinuous points. M. Kobayashi [2] considered the symmetrical situation, indicated that one spectrum could uniquely determine the symmetric potential and the symmetric jump conditions. When there are much more discontinuous points, the situation would become very complicated. M.Shahriari, A.J.Akbarfam, G.Teschl[7] used Weyl function, considered the similar problem with finite discontinuous points.

As to Sturm-Liouville problem on graph, there are more and more researches are under concerning since it is the natural expansion of classical Sturm-Liouville problems and it also relates to many physical problems. Many similar results have been discovered ([3,4,5,6,7].)

\section{\large{Main Result}}

Consider the Sturm-Liouville problem on three-star graph
\begin{gather}
-y''_{j}(x_j)+q_j (x_j)y_j (x_j)=\lambda y_j(x_j), x_j \in (0, \pi) ,j=1,2,3,\\
y'_j(0,\lambda)-h_j y_j (0,\lambda)=0,j=1,2,3,\\
y_1(\pi,\lambda)=y_2(\pi,\lambda)=y_3(\pi,\lambda),\\
y'_1(\pi,\lambda)+y'_2(\pi,\lambda)+y'_3(\pi,\lambda)=0,
\end{gather}
and its jump conditions at $x_1=d$
\begin{gather}
y_1(d+0,\lambda)=ay_1(d-0,\lambda),\\
y'_1(d+0,\lambda)=\frac{1}{a}y_1'(d-0,\lambda)+by_1(d-0,\lambda).
\end{gather}
Denote problem (2.1)-(2.6) by $L=L(q_1,q_2,q_3;h_1,h_2,h_3;a,d,b)$, where $h_j\neq0 ,q_j \in L[0,\pi],j=1,2,3,a>0,|a-1|+|b|>0.$ (This assumption is natural, otherwise it would have no jump at $d$. ) $q_j$ are real, $j=1,2,3.$

Denote the spectrum of the problem
\begin{equation}
\left\{
\begin{array}{lr}
-y''_{i}(x_i)+q_i (x_i)y_i (x_i)=\lambda y_i(x_i), x_i \in (0,\pi), & \\
y'_i(0,\lambda)-h_i y_i (0,\lambda)=0=y_i(\pi,\lambda),&
\end{array}
\right.\nonumber
\end{equation}
by $\sigma_i,i=2,3.$  We agree that if a certain symbol $\delta$ denotes an object relates to $L$, then $\widetilde{\delta}$ would denote an analogous object relates to $ \widetilde{L}=\widetilde{L}(\widetilde{q}_1,\widetilde{q}_2,\widetilde{q}_3;\widetilde{h}_1,\widetilde{h}_2,\widetilde{h}_3;\widetilde{a},\widetilde{d},\widetilde{b})$ .

\begin{theorem} \rm Denote the set of eigenvalues of problem (2.1)-(2.6) by $\Omega$, assume $d<\frac{\pi}{2}, \sigma_i  \cap  \Omega = \emptyset,i=2,3$, $h_i= \widetilde{h}_i,i=2,3,$ on the interval $\ (0,\pi)$ , $q_i(x)= \widetilde{q} _i(x) \ a.e, i=2,3,$ and on the interval $\ (\frac{\pi}{2},\pi)$ , $q_1(x)= \widetilde{q}_1(x) \ a.e.  $.  If $\Omega=  \widetilde{\Omega}$, then on the interval $\ (0,\pi)$ , $q_1(x)=\widetilde{q} _1(x)   \  a.e.  , and \ h_1=\widetilde{h}_1,a=\widetilde{a},b=\widetilde{b},d=\widetilde{d}$.
\end{theorem}
To prove theorem 2.1, we need to study the property of its characteristic function. Denote the solutions of
\begin{equation}
\left\{
\begin{array}{lr}
-y''_{j}(x_j)+q_j (x_j)y_j (x_j)=\lambda y_j(x_j), x_j \in (0,\pi),j=1,2,3, & \\
y'_j(0,\lambda)=h_j ,y_j (0,\lambda)=1,j=1,2,3,&
\end{array}
\right.
\end{equation}
by $\varphi _j(x,\lambda),j=1,2,3.$   Especially, when $j=1$,  $\varphi _1(x,\lambda)$ satisfies the conditions (2.5)-(2.6). Thus, $\lambda$ is an eigenvalue of $ L $ if and only if $\lambda$ is a zero of
\[\omega(\lambda)=
\left |\begin{array}{cccc}
{\varphi}_1(\pi,\lambda) & -\varphi _2 (\pi,\lambda) &0 \\
{\varphi}_1(\pi,\lambda) & 0 & -\varphi _3 (\pi,\lambda)\\
{\varphi}'_1(\pi,\lambda) & \varphi'_2 (\pi,\lambda) &  \varphi'_3(\pi,\lambda)
\end{array}\right |.
\]
  $\omega(\lambda)$ is named the characteristic function of $L$, next, we consider the distribution of its zeros. Let $\lambda =k^2$ ,denote $R_n$ be a rectangle whose vertexes are $\pm(n-\frac{1}{4})+i0$ and $\pm(n-\frac{1}{4})+i(n-\frac{1}{4})$ on $k-plane$ , denote $\Gamma_n$  be the contour on $\lambda$-plane corresponding with $R_n$ and $\Im k>0$. Denote $S_n$ be a rectangle whose vertexes are $\pm(n+\frac{1}{4})+i0$ and $\pm(n+\frac{1}{4})+i(n-\frac{1}{4})$ on $k$-plane, $\gamma_n$ be the contour on $\lambda -plane $ corresponding with $S_n$ and $\Im k>0$ .

\begin{lemma} \rm $\omega(\lambda)$ is an entire function of growth order $\frac{1}{2}$, its zeros are all real, if $\frac{a^2-1}{(a^2+1)}<\frac{3}{4},$ then its zeros could be denoted by $\{\lambda_{n,1},\lambda_{n,2},\lambda_{n,3}\}_{n=1}^\infty,$ and satisfy
\begin{align}
\sqrt{\lambda_{n,1}}=n-1+O(1), \sqrt{\lambda_{n,i}}=n-\frac{1}{2}+O(1), i=2,3.
\end{align}
\end{lemma}
Prove£ºBecause of the property of symmetric operators, it is not difficult to deduce the zeros of $\omega(\lambda)$ are all real. Let $\lambda =k^2, k=\sigma+i\tau$, O.H.Hald [1] proved
\begin{align}
\varphi_i(x,\lambda)= &\cos kx+\frac{h_i}{k}\sin kx+\int_0^x \frac{\sin k(x-t)}{k} q_i(t) \varphi _i(t,\lambda) dt\nonumber\\
=&\cos kx +O(\frac{e^{|\tau| x}}{k}),x\in(0,\pi),i=2,3;\nonumber\\
{\varphi}_1(x,\lambda)=&a \cos k(x-d)\cos kd -\frac{1}{a}\sin k(x-d)\sin kd+\frac{b}{k}\sin k(x-d)\cos kd\nonumber\\
&+\frac{h_1}{k}[\frac{1}{a}\sin k(x-d)\cos kd+a\cos k(x-d)\sin kd+\frac{b}{k}\sin k(x-d)\sin kd]\nonumber\\
&+\frac{1}{k}\int_0^d[\frac{1}{a}\sin k(x-d)\cos k(d-t)+a\cos k(x-d)\sin k(d-t)\nonumber\\
&+\frac{b}{k}\sin k(x-d)\sin k(d-t)]q_1(t)\varphi_1(t,\lambda)dt+\frac{1}{k}\int_d^x\sin k(x-t)q_1(t)\varphi_1(t,\lambda)dt\nonumber\\
=&a \cos k(x-d)\cos kd -\frac{1}{a}\sin k(x-d)\sin kd+O(\frac{e^{\mid\tau\mid x}}{k})\nonumber\\
=&\frac{a+\frac{1}{a}}{2} \cos kx+\frac{a-\frac{1}{a}}{2} \cos k(x-2d)+O(\frac{e^{|\tau| x}}{k})\nonumber\\
=&\alpha_1\cos kx +\alpha_2\cos k(x-2d)+O(\frac{e^{|\tau| x}}{k}),x\in(d,\pi),\nonumber
\end{align}
and
\begin{align}
\varphi'_i(x,\lambda)=&-k\sin kx+O(e^{|\tau| x}),x\in(0,\pi),i=2,3,\nonumber\\
{\varphi}'_1(x,\lambda)=&-(ka\sin k(x-d)\cos kd +\frac{k}{a}\cos k(x-d) \sin kd)+b\cos k(x-d)\cos kd\nonumber\\
&+h_1[\frac{1}{a}\cos k(x-d)\cos kd-a\sin k(x-d)\sin kd+\frac{b}{k}\cos k(x-d)\sin kd]\nonumber\\
&+\int_0^d[\frac{1}{a}\cos k(x-d)\cos k(d-t)-a\sin k(x-d)\sin k(d-t)\nonumber\\
&+\frac{b}{k}\cos k(x-d)\sin k(d-t)]q_1(t)\varphi_1(t,\lambda)dt+\int_d^x \cos k(x-t)q_1(t)\varphi_1(t,\lambda)dt\nonumber\\
=&-(ka\sin k(x-d)\cos kd +\frac{k}{a}\cos k(x-d) \sin kd)+O(e^{|\tau|x})\nonumber\\
=&-k(\alpha_1 \sin kx+\alpha_2 \sin k(x-2d))+O(e^{|\tau |x}),x\in(d,\pi),\nonumber
\end{align}
where $\alpha_1=\frac{a+\frac{1}{a}}{2},\alpha_2=\frac{a-\frac{1}{a}}{2},$ thus
\begin{align}
\omega (\lambda)&=-k[3\alpha_1 \sin k\pi \cos^2 k\pi +\alpha_2\sin k(\pi-2d) \cos^2 k\pi +2\alpha_2 \cos k(\pi-2d)\sin k\pi \cos k\pi]+O(e^{3|\tau| \pi})\nonumber\\
&=-\frac{3}{2}\alpha_1 k \cos k \pi[\sin 2k\pi+3\beta \sin2k(\pi-d)+\beta\sin2kd]+O(e^{3|\tau |\pi})\nonumber\\
&=\omega_0(\lambda)+O(e^{3|\tau |\pi}),\nonumber
\end{align}
where $\beta=\frac{\alpha_2}{3\alpha_1}=\frac{a^2-1}{3(a^2+1)},\omega_0(\lambda)=-\frac{3}{2}\alpha_1 k \cos k \pi[\sin 2k\pi+3\beta \sin2k(\pi-d)+\beta\sin2kd]$.
Consider $k=\sigma +i\tau=\pm(n-\frac{1}{4})+i\tau$ first,
it is obvious that there exists a positive number $c$ s.t. $$|\cos k\pi|>ce^{\tau\pi}.$$
Since
\begin{align}
|\sin 2k\pi+3\beta \sin2k(\pi-d)+\beta\sin2kd|&\geq \cosh2\tau\pi-4\beta\cosh2\tau\pi\geq \frac{1-4\beta}{2}e^{2\tau\pi},\nonumber
\end{align}
and because of the assumption, we have $\beta < \frac{1}{4}$, so
\begin{align}
|\omega_0(\lambda)|\geq c|k|e^{3\tau \pi}.
\end{align}

Then consider when $k=\sigma +i\tau=(n-\frac{1}{4})i+\sigma,$ we still have $|\cos k\pi|>ce^{\tau\pi}. $
Similarly , it is not difficult to prove the following when $n$ is large enough,
\begin{align}
|\sin 2k\pi+3\beta \sin2k(\pi-d)+\beta\sin2kd|\geq ce^{2\tau\pi},\nonumber
\end{align}
Thus we still have (2.9) . Similarly, we could deduce (2.9) when $k \in S_n$ , so by $ Rouch\acute{e} $ Theorem, $\omega(\lambda)$  and $\omega_0(\lambda)$ have the same number of zeros in $\Gamma_n$ and $\gamma_n$ when $n$  is large enough. Thus, to find out the distribution of zeros of $\omega(\lambda)$, we need to study the distribution of zeros of $\omega_0(\lambda)$ at first.

From the discussion above, we have
\begin{align}
|\sin 2k\pi|>|3\beta \sin2k(\pi-d)+\beta\sin2kd|,\nonumber
\end{align}
when $k \in R_n$  or $k \in S_n$ ($n$ is large enough), thus $\omega_0(\lambda)$ and $z(\lambda)=k\cos k\pi \sin2k\pi$ have the same number of zeros in $\Gamma_n$, and $\gamma_n$ . And it is obvious that $z(\lambda)$ has $3n$ and $3n+1$ zeros in $\Gamma_n$ and $\gamma_n$, respectively. So we could deduce that $\omega(\lambda)$ has $3n$ and $3n+1$ zeros in $\Gamma_n$ and $\gamma_n$, respectively. And we have indicated that its zeros are all real at the beginning, thus we get the left half of (2.8). For the right half, we need to repeat the above steps on the rectangle whose vertexes are $\pm(n+1-\frac{1}{4})+i0$ and $\pm(n+1-\frac{1}{4})+i(n-\frac{1}{4})$ on $k-plane$ and the corresponding contour on $\lambda-plane $, we omit the steps.

\begin{lemma}\rm Let
\begin{align}
\varphi_0(x,\lambda)&=\cos kx ,&0\leq x<d,\nonumber\\
&=a\cos k(x-d)\cos kd-\frac{1}{a}\sin k(x-d)\sin kd,&d<x \leq 2\pi.\nonumber
\end{align}
then there exists a bounded function $K_j(x,t), j=1,2,3,$ such that $$\varphi_1(x,\lambda)=\varphi_0(x,\lambda)+\int_0^x K_1(x,t)\cos kt dt,$$
$$\varphi_i(x,\lambda)=\cos kx+\int_0^x K_i(x,t)\cos kt dt,i=2,3.$$
Let $K_j(x,t)=0$ if $ x<t$ or $t<0, j=1,2,3$, $K_2(x,t), K_3(x,t) $ are continuous and bounded on the area $0\leq x \leq 2\pi,0\leq t \leq x$ , $K_1(x,t)$ is bounded and continuous on the areas $K_1, K_2, K_3, K_4, K_5$ , where
\begin{align}
&K_1: 0\leq x \leq d,  0\leq t \leq x,\nonumber\\
&K_2: d< x <2 d, 0\leq t \leq 2d-x,\nonumber\\
&K_3: d< x <2 d, 2d-x< t \leq x,\nonumber\\
&K_4: 2d\leq x \leq\pi,x-2d\leq t \leq x,\nonumber\\
&K_5: 2d< x \leq\pi, 0\leq t < x-2d.\nonumber
\end{align}
\end{lemma}

Proof was given in [1].

The following lemma uses a number of trigonometric function transformations, here we deduce an easy one, and the others could be deduced by the same way. Fix the positive number $A$ large enough,
\begin{align}
&\frac{1}{T^2}\int_A^T x\cos^2 ax\sin b_1 x \sin b_2x dx\nonumber\\
=&\frac{1}{T^2}\int_A^T x\frac{1+\cos 2ax}{2}\frac{\cos (b_1-b_2)x-\cos(b_1+b_2)x}{2}dx\nonumber\\
=&\frac{1}{4T^2}\int_A^T[x \cos(b_1-b_2)x-x \cos(b_1+b_2)x+x\cos2ax\cos(b_1-b_2)x-x\cos2ax\cos(b_1+b_2)x]dx\nonumber\\
=&\frac{1}{4T^2}\int_A^T\{x \cos(b_1-b_2)x-x \cos(b_1+b_2)x+\frac{x}{2}[\cos(2a+(b_1-b_2))x\nonumber\\
&+\cos(2a-(b_1-b_2))x]-\frac{x}{2}[\cos(2a+(b_1+b_2))x+\cos(2a-(b_1+b_2))x]\}dx\nonumber
\end{align}
for the above four terms, it is not difficult to deduce that when $T$ tends to infinity,
\begin{equation}
\frac{1}{4T^2}\int_A^Tx \cos(b_1-b_2)xdx= \left\{
\begin{array}{lr}
\frac{1}{8}+O(\frac{1}{T}),
&b_1-b_2=0,\\
O(\frac{1}{T}),&b_1-b_2\neq0,
\end{array}
\right.\nonumber
\end{equation}
\begin{equation}
\frac{1}{4T^2}\int_A^Tx \cos(b_1+b_2)xdx= \left\{
\begin{array}{lr}
\frac{1}{8}+O(\frac{1}{T}),
&b_1+b_2=0,\\
O(\frac{1}{T}),&b_1+b_2\neq0,
\end{array}
\right.\nonumber
\end{equation}
\begin{equation}
\frac{1}{4T^2}\int_A^T \frac{x}{2}[\cos(2a+(b_1-b_2))x+\cos(2a-(b_1-b_2))x]dx= \left\{
\begin{array}{lr}
\frac{1}{16}+O(\frac{1}{T}),
&(2a)^2-(b_1-b_2)^2=0,\\
O(\frac{1}{T}),&(2a)^2-(b_1-b_2)^2\neq0,
\end{array}
\right.\nonumber
\end{equation}
\begin{equation}
\frac{1}{4T^2}\int_A^T \frac{x}{2}[\cos(2a+(b_1+b_2))x+\cos(2a-(b_1+b_2))x]dx= \left\{
\begin{array}{lr}
\frac{1}{16}+O(\frac{1}{T}),
&(2a)^2-(b_1+b_2)^2=0,\\
O(\frac{1}{T}),&(2a)^2-(b_1+b_2)^2\neq0,
\end{array}
\right.\nonumber
\end{equation}

\begin{lemma} \rm If $\frac{a^2-1}{(a^2+1)}<\frac{3}{4}, |a-1|+|b|>0,d<\frac{\pi}{2},$ then $a$ and $ d$ in L could be decided by $\Omega$.  Moreover, if $a\neq 1$, then we do not need the assumption $d<\frac{\pi}{2}$ .
\end{lemma}
Proof. From Lemma2.2 and Hadamard factorization theorem, there exists $C$ such that $\omega(\lambda)=C\widetilde{\omega}(\lambda)$, thus $\omega_0(\lambda)-C\widetilde{\omega}_0(\lambda)=C(\widetilde{\omega}(\lambda)-\widetilde{\omega}_0(\lambda))-(\omega(\lambda)-\omega_0(\lambda))$.
So when $k$ is a real number, we have
\begin{align}
&C\widetilde{\alpha}_1 k\cos k\pi \sin 2k\pi -\alpha_1 k\cos k\pi\sin2k\pi\nonumber\\
+&3C\widetilde{\alpha}_1 \widetilde{\beta}k\cos k\pi \sin 2k(\pi-\widetilde{d})-3\alpha_1 \beta k\cos k\pi\sin2k(\pi-d)\nonumber\\
+&C\widetilde{\alpha}_1\widetilde{\beta}k\cos k\pi \sin2k\widetilde{d}-\alpha_1 \beta k\cos k\pi\sin 2kd=O(1).
\end{align}
First we conduct the following operation on (2.10)$$\frac{1}{T^2}\int_A^T *\cos k\pi\sin2k\pi dk$$ we get ($*$ represents each term of (2.10) )
$$C\widetilde{\alpha}_1(\frac{1}{8}+O(\frac{1}{T}))-\alpha_1(\frac{1}{8}+O(\frac{1}{T}))+O(\frac{1}{T})=O(\frac{1}{T}),$$
so we have $C\widetilde{\alpha}_1=\alpha_1$, and then (2.10) changes to
\begin{align}
&3 \widetilde{\beta}k\cos k\pi \sin 2k(\pi-\widetilde{d})-3 \beta k\cos k\pi\sin2k(\pi-d)\nonumber\\
&+\widetilde{\beta}k\cos k\pi \sin2k\widetilde{d}- \beta k\cos k\pi\sin 2kd=O(1).
\end{align}

We divide the following discussion into two parts, $d=\widetilde{d}$ and $d\neq \widetilde{d}$.

I1. $d=\widetilde{d}=\frac{\pi}{2}.$ We conduct the following operation on (2.11) $$\frac{1}{T^2}\int_A^T *\cos k\pi\sin k\pi dk$$ then we could get
$$(4 \widetilde{\beta}-4\beta)(\frac{1}{8}-\frac{1}{16}+O(\frac{1}{T}))=O(\frac{1}{T}),$$
so $\widetilde{\beta}=\beta,$ and then we find $\widetilde{a}=a, C=1.$

I2. $d=\widetilde{d}\neq \frac{\pi}{2}.$ Similarly, we conduct the operation $\frac{1}{T^2}\int_A^T *\cos k\pi\sin 2k(\pi-d) dk$ on (2.11) could get
$$(3\widetilde{\beta}-3\beta)(\frac{1}{8}+O(\frac{1}{T}))+(\widetilde{\beta}-\beta)(-\frac{1}{16}+O(\frac{1}{T}))=O(\frac{1}{T}),$$
so we also have $\widetilde{\beta}=\beta,$ $\widetilde{a}=a, C=1.$

Thus, for the part $d=\widetilde{d}$, we proved the lemma. Then we consider the second part which is much more complicated.

II1. $\widetilde{d}=\frac{\pi}{2}, d\neq \frac{\pi}{2}.$ We conduct the following operation on (2.11)$$\frac{1}{T^2}\int_A^T *\cos k\pi\sin k\pi dk$$
get
$$4 \widetilde{\beta}(\frac{1}{8}-\frac{1}{16}+O(\frac{1}{T}))-O(\frac{1}{T})=O(\frac{1}{T}),$$
thus $\widetilde{\beta}=0,$ then $\widetilde{a}=1$. Then repeat the similar steps in I2, we could get $\beta=0,a=1$.

II2. $\widetilde{d},d \neq \frac{\pi}{2}, d+\widetilde{d} \neq \pi.$ First, we conduct the  operation $\frac{1}{T^2}\int_A^T *\cos k\pi\sin 2k(\pi-\widetilde{d}) dk$ on (2.11) could get
$$3\widetilde{\beta}(\frac{1}{8}+O(\frac{1}{T}))+\widetilde{\beta}(-\frac{1}{16}+O(\frac{1}{T})=O(\frac{1}{T}),$$
Thus we have $\widetilde{\beta}=0,\widetilde{a}=1$. Then we conduct the  operation $\frac{1}{T^2}\int_A^T *\cos k\pi\sin 2k(\pi-d) dk$ on (2.11), then we have $\beta=0, a=1$.

II3. $\widetilde{d},d \neq \frac{\pi}{2}, d+\widetilde{d} = \pi.$ Similar to II2, it is not difficult to get $a=\widetilde{a}=1$.

We are going to prove that we would have $b=0$ if $d\neq \widetilde{d}$£¬then $|a-1|+|b|=0$, thus it is contradictory, so $d=\widetilde{d}$.  During the following, we assume $d, \widetilde{d}<\frac{\pi}{2}$. We study more about the characteristic function, since $$\omega(\lambda)=\sum_{j=1}^{3}\varphi'_j (\pi,\lambda)\prod_{i=1,2,3,i\neq j}\varphi_i(\pi,\lambda):=\sum_{j=1}^{3}F_j(\lambda),$$
 where $$F_j(\lambda)=\varphi'_j (\pi,\lambda)\prod_{i=1,2,3,i\neq j}\varphi_i(\pi,\lambda),j=1,2,3.$$
since $a=1$, from Lemma 2.2, Lemma2.3, we have
\begin{align}
F_1(\lambda)=&\varphi'_1(\pi,\lambda)\varphi_2(\pi,\lambda)\varphi_3(\pi,\lambda)\nonumber\\
=&[-k\sin k\pi+(\frac{b}{2}+h_1+\frac{1}{2}\int_0^{\pi}q_1(t)dt)\cos k\pi+\frac{b}{2}\cos k(\pi-2d)\nonumber\\&+\frac{1}{2}\int_0^{\pi}\cos k(\pi-2t)q_1(t)dt+E_1(k)]\nonumber\\
&(\cos k\pi +\int_0^{\pi}K_2(t)\cos ktdt)(\cos k\pi +\int_0^{\pi}K_3(t)\cos ktdt),\nonumber
\end{align}
where
\begin{align}
E_1(k)=&\frac{h_1b}{k}\cos k(\pi-d)\sin kd+\int_0^{\pi}\cos k(\pi-t)q_1(t)(\varphi_1(t,\lambda)-\cos kt)dt\nonumber\\&+\frac{b}{k}\int_0^d\cos k(\pi-d)\sin k(d-t)q_1(t)\varphi_1(t,\lambda)dt,\nonumber
\end{align}
by Lemma2.2, and Paley-Wiener Theorem, similar to [1, Lemma6], we have
\begin{align}
F_1(\lambda)=&[-k\sin k\pi+A_1\cos k\pi+\frac{b}{2}\cos k(\pi-2d)+\int_0^{\pi}P_1(t)\cos ktdt]\nonumber\\
&(\cos^2 k\pi +\cos k\pi \int_0^{\pi}Z_1(t)\cos ktdt+\int_0^{2\pi}Z_2(t)\cos ktdt),\nonumber
\end{align}
where $A_1=\frac{b}{2}+h_1+\frac{1}{2}\int_0^{\pi}q_1(t)dt$, $P_1(t)$ derives from $E_1(k)$ and $\frac{1}{2}\int_0^{\pi}\cos k(\pi-2t)q_1(t)dt$, $Z_1(t)=K_2(t)+K_3(t)$, also,
\begin{align}
&\int_0^{\pi}K_2(t)\cos ktdt\int_0^{\pi}K_3(t)\cos ktdt=\int_0^{\pi}\int_0^{\pi}K_2(t)K_3(s)\cos kt\cos ksdtds\nonumber\\
=&\frac{1}{2}\int_0^{\pi}\int_s^{s+\pi}K_2(x-s)K_3(s)\cos kx dxds+\frac{1}{2}\int_0^{\pi}\int_{-s}^{\pi-s}K_2(x+s)K_3(s)\cos kxdxds\nonumber\\
=&\frac{1}{2}\int_0^{\pi}\int_0^xK_2(x-s)K_3(s)ds\cos kxdx+\frac{1}{2}\int_{\pi}^{2\pi}\int_{x-\pi}^{\pi}K_2(x-s)K_3(s)ds\cos kxdx\nonumber\\
&+\frac{1}{2}\int_0^{\pi}\int_0^{\pi-x}K_2(x+s)K_3(s)ds\cos kxdx+\frac{1}{2}\int_0^{\pi}\int_x^{\pi}K_2(-x+s)K_3(s)ds\cos kx dx,\nonumber
\end{align}
thus, we could get $Z_2(t)$ by turning the letter $x$ into $t$ . Similarly,
\begin{align}
F_2(\lambda)=&(\cos k\pi+\int_0^{\pi}K_1(t)\cos ktdt)(\cos k\pi+\int_0^{\pi}K_3(t)\cos ktdt)\nonumber\\ &[-k\sin k\pi+(h_2+\frac{1}{2}\int_0^{\pi}q_2(t)dt)\cos k\pi+\frac{1}{2}\int_0^\pi\cos k(\pi-2t)q_2(t)dt\nonumber\\&+\int_0^\pi\cos k(\pi-t)q_2(t)(\varphi_2(t,\lambda)-\cos kt)dt]\nonumber
\nonumber\\=&(-k\sin k\pi+A_2\cos k\pi+\int_0^{\pi}P_2(t)\cos ktdt)[\cos^2k\pi+\cos k\pi\int_0^{\pi}Z_3(t)\cos ktdt\nonumber\\&+\int_0^{2\pi}Z_4(t)\cos ktdt],\nonumber
\end{align}
\begin{align}
F_3(\lambda)=&(\cos k\pi+\int_0^{\pi}K_1(t)\cos ktdt)(\cos k\pi+\int_0^{\pi}K_2(t)\cos ktdt)\nonumber\\ &[-k\sin k\pi+(h_3+\frac{1}{2}\int_0^{\pi}q_3(t)dt)\cos k\pi+\frac{1}{2}\int_0^\pi\cos k(\pi-2t)q_3(t)dt\nonumber\\&+\int_0^\pi\cos k(\pi-t)q_3(t)(\varphi_3(t,\lambda)-\cos kt)dt]\nonumber
\nonumber\\=&(-k\sin k\pi+A_3\cos k\pi+\int_0^{\pi}P_3(t)\cos ktdt)[\cos^2k\pi+\cos k\pi\int_0^{\pi}Z_5(t)\cos ktdt\nonumber\\&+\int_0^{2\pi}Z_6(t)\cos ktdt],\nonumber
\end{align}
where $A_i=h_i+\frac{1}{2}\int_0^{\pi}q_i(t)dt $; $P_i(t), Z_{2i-1}(t),$ and $ Z_{2i}(t) $ derive from the similar way like $P_1(t), Z_1(t),Z_2(t)$ ,respectively, i=2,3.
Thus, we have
\begin{align}
\omega(\lambda)=&-3k\sin k\pi\cos^2k\pi+(A_1+A_2+A_3)\cos^3k\pi\nonumber\\&+\frac{b}{2}\cos k(\pi-2d)(\cos^2 k\pi +\cos k\pi \int_0^{\pi}Z_1(t)\cos ktdt+\int_0^{2\pi}Z_2(t)\cos ktdt))\nonumber\\
&-k\sin k\pi \cos k\pi\int_0^{\pi}(Z_1(t)+Z_3(t)+Z_5(t))\cos ktdt\nonumber\\
&-k\sin k\pi\int_0^{2\pi}(Z_2(t)+Z_4(t)+Z_6(t))\cos ktdt\nonumber\\
&+\cos k\pi \int_0^{2\pi}(A_1Z_2(t)+A_2Z_4(t)+A_3Z_6(t))\cos ktdt\nonumber\\
&+\cos^2 k\pi \int_0^{\pi}(A_1Z_1(t)+A_2Z_3(t)+A_3Z_5(t))\cos ktdt\nonumber\\
&+\cos^2 k\pi\int_0^{\pi}(P_1(t)+P_2(t)+P_3(t))\cos ktdt\nonumber\\
&+\cos k\pi[\int_0^{\pi}P_1(t)\cos ktdt\int_0^{\pi}Z_1(t)\cos ktdt+\int_0^{\pi}P_2(t)\cos ktdt\int_0^{\pi}Z_3(t)\cos ktdt\nonumber\\
&+\int_0^{\pi}P_3(t)\cos ktdt\int_0^{\pi}Z_5(t)\cos ktdt]\nonumber\\
&+\int_0^{\pi}P_1(t)\cos ktdt\int_0^{2\pi}Z_2(t)\cos ktdt+\int_0^{\pi}P_2(t)\cos ktdt\int_0^{2\pi}Z_4(t)\cos ktdt\nonumber\\
&+\int_0^{\pi}P_3(t)\cos ktdt\int_0^{2\pi}Z_6(t)\cos ktdt\\
=&-3k\sin k\pi\cos^2k\pi+B\cos^3k\pi\nonumber\\
&+\frac{b}{2}\cos k(\pi-2d)(\cos^2 k\pi +\cos k\pi \int_0^{\pi}Z_1(t)\cos ktdt+\int_0^{2\pi}Z_2(t)\cos ktdt))\nonumber\\
&-k\sin k\pi \cos k\pi\int_0^{\pi}L_1(t)\cos ktdt
-k\sin k\pi\int_0^{2\pi}L_2(t)\cos ktdt\nonumber\\
&+\cos k\pi \int_0^{2\pi}L_3(t)\cos ktdt+\cos^2 k\pi \int_0^{\pi}L_4(t)\cos ktdt+\int_0^{3\pi}L_5(t)\cos ktdt,
\end{align}
where $B=A_1+A_2+A_3,L_1(t)=Z_1(t)+Z_3(t)+Z_5(t),L_2(t)=Z_2(t)+Z_4(t)+Z_6(t),$ $L_4(t)=A_1Z_1(t)+A_2Z_3(t)+A_3Z_5(t)+P_1(t)+P_2(t)+P_3(t)$,  and $L_3(t)$ derives from the fifth, eighth, and ninth row of (2.12) , respectively. $L_5(t)$ derives from the last two rows of (2.12). Also, $K_2(t),K_3(t)$ could be represented as the sum of an absolutely continuous function and a function which is continuously differentiable with respect to $t$ . ( See the proof of [1, Lemma5] ). So $Z'_2(t)\in L[0,2\pi]$, similarly, $Z'_4(t), Z'_6(t)\in L[0,2\pi]$, thus $L'_2(t)\in L[0,2\pi]$.

By (2.13),
\begin{align}
0=&\omega(\lambda)-\widetilde{\omega}(\lambda)\nonumber\\=&(B-\widetilde{B})\cos^3k\pi+\frac{b}{2}\cos k(\pi-2d)\cos^2 k\pi-\frac{\widetilde{b}}{2}\cos k(\pi-2\widetilde{d})\cos^2 k\pi\nonumber\\
&+\frac{b}{2}\cos k(\pi-2d)\cos k\pi\int_0^{\pi}Z_1(t)\cos ktdt-\frac{\widetilde{b}}{2}\cos k(\pi-2\widetilde{d})\cos k\pi\int_0^{\pi}\widetilde{Z}_1(t)\cos ktdt\nonumber\\
&+\frac{b}{2}\cos k(\pi-2d)\int_0^{2\pi}Z_2(t)\cos ktdt-\frac{\widetilde{b}}{2}\cos k(\pi-2\widetilde{d})\int_0^{2\pi}\widetilde{Z}_2(t)\cos ktdt\nonumber\\
&-k\sin k\pi \cos k\pi\int_0^{\pi}(L_1(t)-\widetilde{L}_1(t))\cos ktdt-k\sin k\pi\int_0^{2\pi}(L_2(t)-\widetilde{L}_2(t))\cos ktdt\nonumber\\
&+\cos k\pi \int_0^{2\pi}(L_3(t)-\widetilde{L}_3(t))\cos ktdt+\cos^2 k\pi \int_0^{\pi}(L_4(t)-\widetilde{L}_4(t))\cos ktdt\nonumber\\
&+\int_0^{3\pi}(L_5(t)-\widetilde{L}_5(t))\cos ktdt,
\end{align}
 conduct the following operation on (2.14)
$$\frac{1}{T}\int_A^T*\cos k(\pi-2d)\cos^2 k\pi dk=\frac{1}{T}\int_A^T *\xi dk,$$
where $\xi$ represents $\cos k(\pi-2d)\cos^2 k\pi,$ by Fubini Theorem, it is not difficult to deduce the following,
\begin{align}
&\frac{1}{T}\int_A^T\cos^3k\pi \xi dk=O(\frac{1}{T}),\frac{1}{T}\int_A^T\cos k(\pi-2d)\cos^2 k\pi \xi dk=\frac{3}{16}+O(\frac{1}{T}),\nonumber\\
&\frac{1}{T}\int_A^T\cos k(\pi-2\widetilde{d})\cos^2 k\pi \xi dk=O(\frac{1}{T}),\nonumber\\
&\frac{1}{T}\int_A^T\cos k(\pi-2d)\cos k\pi\int_0^{\pi}Z_1(t)\cos ktdt\xi dk=O(\frac{1}{T}),\nonumber\\
&\frac{1}{T}\int_A^T\cos k(\pi-2\widetilde{d})\cos k\pi\int_0^{\pi}\widetilde{Z}_1(t)\cos ktdt\xi dk=O(\frac{1}{T}),\nonumber\\
&\frac{1}{T}\int_A^T\cos k(\pi-2d)\int_0^{2\pi}Z_2(t)\cos ktdt\xi dk=O(\frac{1}{T}),\nonumber\\&\frac{1}{T}\int_A^T\cos k(\pi-2\widetilde{d})\int_0^{2\pi}\widetilde{Z}_2(t)\cos ktdt\xi dk=O(\frac{1}{T}),\nonumber\\
&\frac{1}{T}\int_A^Tk\sin k\pi \cos k\pi\int_0^{\pi}(L_1(t)-\widetilde{L}_1(t))\cos ktdt \xi dk=O(\frac{1}{T}),\nonumber\\
&\frac{1}{T}\int_A^Tk\sin k\pi\int_0^{2\pi}(L_2(t)-\widetilde{L}_2(t))\cos ktdt\xi dk=O(\frac{1}{T}),\nonumber\\
&\frac{1}{T}\int_A^T\{\cos k\pi \int_0^{2\pi}(L_3(t)-\widetilde{L}_3(t))\cos ktdt+\cos^2 k\pi \int_0^{\pi}(L_4(t)-\widetilde{L}_4(t))\cos ktdt\nonumber\\
&+\int_0^{3\pi}(L_5(t)-\widetilde{L}_5(t))\cos ktdt\}\xi dk=O(\frac{1}{T}).\nonumber
\end{align}
Among the above equations , the equation in the second row needs the assumption $d<\frac{\pi}{2}$ , otherwise if $d+\widetilde{d}=\pi$, then $$\frac{1}{T}\int_A^T\cos k(\pi-2\widetilde{d})\cos^2 k\pi \xi dk=\frac{3}{16}+O(\frac{1}{T}),$$
thus we could not deduce $b=0$ to get a contradiction. The penultimate equation needs integration by parts followed by using Fubini Theorem, its feasibility has been discussed between (2.13) and (2.14) . The others could be deduced by trigonometric function transformation.

Thus, $b=0$£¬ similarly we have $\widetilde{b}=0$, this is contradictory.

\begin{lemma} \rm Let $\varphi_{1-}=\varphi_1(d-0,\lambda), \widetilde{\varphi}_{1-}=\widetilde{\varphi}_1(\widetilde{d}-0,\lambda)$. If $\Omega=\widetilde{\Omega},$ then $a=\widetilde{a},d=\widetilde{d}$, assume $\Omega\cap \sigma_i= \emptyset ,i=2,3,$ and on the interval $(\frac{\pi}{2},\pi)$, $q_1(x)=\widetilde{q}_1(x) \ a.e. , $$h_2=\widetilde{h}_2, h_3=\widetilde{h}_3$, then we have
\begin{gather}
b-\widetilde{b}=-\frac{a^2-a^{-2}}{2a}\int_d^{\frac{\pi}{2}}(q_1-\widetilde{q}_1)dt,\nonumber\\
h_1-\widetilde{h}_1=-\frac{1}{2}\int_0^d (q_1-\widetilde{q}_1)dt-\frac{1}{2a^2}\int_d^{\frac{\pi}{2}} (q_1-\widetilde{q}_1)dt,\nonumber\\
\int_0^d (q_1-\widetilde{q}_1)(\varphi_1 \widetilde{\varphi}_1-\frac{1}{2})dt+\int_d^{\frac{\pi}{2}} (q_1-\widetilde{q}_1)((\varphi_1 \widetilde{\varphi}_1-\frac{1}{2a^2}-\frac{a^2-a^{-2}}{2}\varphi_{1-} \widetilde{\varphi}_{1-})dt=0.\nonumber
\end{gather}
\end{lemma}
\begin{remark}\rm
Lemma 2.5 is similar to [1, Lemma7], yet the assumption $\Omega\cap \sigma_i= \emptyset ,(i=2,3)$ in the statement is significant, this assumption comes from [9, Theorem2.2], and we will see its significance in the proof.
\end{remark}

Proof. We could get the following immediately from (2.1)-(2.6)
\begin{align}
(h_1-\widetilde{h}_1)+a(b-\widetilde{b})\varphi_{1-}\widetilde{\varphi}_{1-}+\int_0^{\frac{\pi}{2}}(q_1-\widetilde{q}_1)\varphi_1\widetilde{\varphi}_1dt=\varphi'_1(\pi,\lambda)\widetilde{\varphi}_1(\pi,\lambda)-\varphi_1(\pi,\lambda)\widetilde{\varphi}'_1(\pi,\lambda),
\end{align}
Denote the left side of (2.15) by $\Phi(\lambda)$, in [1], $\Phi(\lambda_m)\equiv0$, $\lambda_m$ are the eigenvalues of (2.1)-(2.6) (For simplicity , we denote the three classes of eigenvalues in Lemma2.2 by $\lambda_m$ ). In this lemma, however, We could not have such result without the assumption.

Indeed, let $u(x,\lambda)=(u_1(x,\lambda),u_2(x,\lambda),u_3(x,\lambda))$ be a solution of (2.1)-(2.6) , then there exists $C(\lambda)=(C_1(\lambda),C_2(\lambda),C_3(\lambda))\neq 0$ and independent with $x$ , such that $$u_j(x,\lambda)=C_j(\lambda)\varphi_j(x,\lambda),j=1,2,3,$$
By the assumption , $\varphi_i(\pi,\lambda_m)\neq0, i=2,3,\forall m\in N.$ Thus, by (2.3) and (2.4), $\forall m\in N$, we have
\begin{equation}
\left\{
\begin{array}{lr}
C_1(\lambda_m)\varphi_1(\pi,\lambda_m)= C_2(\lambda_m)\varphi_2(\pi,\lambda_m)=C_3(\lambda_m)\varphi_3(\pi,\lambda_m),& \\
C_1(\lambda_m)\varphi'_1(\pi,\lambda_m)+C_2(\lambda_m)\varphi'_2(\pi,\lambda_m)+C_3(\lambda_m)\varphi'_3(\pi,\lambda_m)=0.&
\end{array}
\right.\nonumber
\end{equation}
similar to the proof of [9, Theorem 2.2] , $\forall m\in N$, $$\frac{\varphi'_1(\pi,\lambda_m)}{\varphi_1(\pi,\lambda_m)}=\frac{\widetilde{\varphi}'_1(\pi,\lambda_m)}{\widetilde{\varphi}_1(\pi,\lambda_m)},$$
thus by (2.15) we have $\Phi(\lambda_m)\equiv0$. The following parts see the proof of [1, Lemma 7].

The proof of \textbf{Theorem 2.1}: By Lemma2.4, Lemma 2.5, [1, Lemma7], the proof of [1, Theorem 1], we could immediately prove Theorem 2.1.



\end{document}